\documentclass[11pt]{article}
\usepackage[dvips]{graphicx}
\usepackage[T2A]{fontenc}
\usepackage[cp1251]{inputenc}
\usepackage[english]{babel}
\usepackage{amssymb}
\usepackage{amsmath}
\usepackage{amscd}
\usepackage{amsthm}
\usepackage{comment}

\begin{document}

\newtheorem{te}{Theorem}[section]
\newtheorem{lem}{Lemma}[section]
\newtheorem{cor}{Corollary}
\newtheorem{df}{Definition}[section]
\newtheorem{exam}{Example}[section]

\begin{center}
{\Large\bf Parshin Residues via Coboundary Operators.}

{\footnotesize by Mikhail Mazin, Stony Brook University.}
\end{center}

\begin{abstract}
The article consists of two main parts: an analog of the theory of Leray coboundary operators for
stratified spaces, and its application to the theory of Parshin residues. The first part is
independent from the second. It uses the theory of Whitney stratifications. The second part is an
application of the first. In particular, a geometric proof of the Parshin's Reciprocity Law for
residues in complex case is given.
\end{abstract}

\section{Introduction.}

Let $X$ be a compact complex curve and $\omega$ be a meromorphic $1$-form on $X.$ In an open
neighborhood of each point $x\in X$ we can write

$$
\omega=f(t)dt, \ f(t) =\sum\limits_{i>N} \lambda_i t^i,
$$
\noindent where $t$ is a local normalizing parameter at $x.$ The coefficient $\lambda_{-1}$ of the
above series does not depend on the choice of parameter $t$ and is called the residue of $\omega$
at $x.$ The residue is non-zero only at finitely many points $\Sigma\subset X,$ where $\omega$ has a
pole. The well-known residue formula says that the sum of residues of $\omega$ over all points of
$\Sigma$ is zero:

$$
\sum\limits_{x\in\Sigma}res_x\omega = 0
$$

Indeed, the residue at $x \in\Sigma$ is equal to the integral of $\omega$ over any sufficiently
small cycle enclosing $x,$ divided by $2\pi i.$ In the complement $X\backslash\Sigma$ the form
$\omega$ is closed, and the sum of cycles is homologous to zero. Thus, the residue formula follows
from the Stokes' Theorem.

Although this proof is topological, the residue itself can be defined purely algebraically and one
can give an algebraic proof of the residue formula which works in a much more general situation, not
only in the case of complex curves (see, for example, \cite{S}, \cite{T}).

In the late $70$'s A. Parshin introduced his notion of multidimensional residue for a rational
$n$-form $\omega$ on an $n$-dimensional algebraic variety $V_n.$ (In \cite{P1} Parshin mostly deals
with the two-dimensional case, then A. Beilinson and V. Lomadze in \cite{B} and \cite{L}
generalized Parshin's ideas to the multidimensional case). The main difference between the Parshin
residue and the classical one-dimensional residue is that in higher dimensions one computes the
residue not at a point, but at a complete flag of subvarieties $F=\{V_n\supset\dots\supset V_0\},$
$\dim V_k=k.$

Parshin, Beilinson, and Lomadze proved the Reciprocity Law for multidimensional residues, which
generalizes the classical residue formula:

Fix a partial flag of irreducible subvarieties
$\{V_n\supset\dots\supset\widehat{V_k}\supset\dots\supset V_0\}$, where $V_k$ is omitted ($0<k<n$).
Then
$$
\sum\limits_{V_{k+1}\supset X\supset V_{k-1}} res_{V_n\supset\dots\supset X\supset\dots\supset
V_0}(\omega)=0,
$$
\noindent where the sum is taken over all irreducible $k$-dimensional subvarieties $X,$ such that
$V_{k+1}\supset X\supset V_{k-1}.$ More precisely, the statement of
the theorem is that there are only finitely many summands which are not zeros, and the sum of
non-zero summands is zero.

In addition, if $V_1$ is proper (compact in the complex case), then one has the same relation for
$k=0:$
$$
\sum\limits_{x\in V_1} res_{V_n\supset\dots\supset V_1\supset \{x\}}(\omega)=0.
$$
\noindent Again, there are finitely many non-zero summands and their sum is zero.

All these papers are purely algebraic. The methods used by Parshin, Beilinson and Lomadze are
applicable in very general settings, not only over complex numbers. However, in the complex case
one would expect a more geometric variant of the theory.

J.-L. Brylinski and D. A. McLaughlin give a more topological treatment of the complex case in
\cite{BrM}. Given a flag $F=\{V_n\supset\dots\supset V_0\}$ they introduce {\it flag-localized}
homology groups $H_*^{V_i}(V_n;F)$ and a homology class $k_F\in H_n^{V_n}(V_n;F),$ such that
$$
res_F \omega=\frac{1}{(2\pi i)^n}\int\limits_{k_F}\omega
$$
\noindent for any meromorphic $n$-form $\omega.$ The class $k_F$ is obtained from of the
fundamental class $c_{V_0}\in H_{2n}(V_n,V_n\backslash V_0)$ by applying the boundary homomorphisms
in the appropriate flag-localized homology groups $n$ times. J.-L. Brylinski and D. A. McLaughlin
mention, that the class $k_F$ could be constructed in a more geometric way, so that it is naturally
represented by a union of certain real $n$-tori. However, they only give such a construction in the
case when all elements of the flag $F$ are smooth.

In this paper we develop a different approach to the construction of the class $k_F.$ We use the
geometry of the Whitney stratified spaces to introduce the Leray coboundary operators
$\phi_{X,Y}:H_*(X)\to H_{*+k-n-1}(Y)$ for any two consecutive strata $X<Y,\ \dim X=n,\  \dim Y=k$ (see Definition \ref{order on strata}) of
a stratified space. Given a flag $F=\{V_n\supset\dots\supset V_0\}$ and a meromorphic top-form
$\omega$ on $V_n,$ one can choose a stratification of $V_n,$ such that the flag $F$ consists of
closures of strata, and $\omega$ is regular on the top-dimensional stratum.

Then one can construct the homology class $\Delta_F:=\phi_{\breve V_{n-1},\breve
V_n}\circ\dots\circ\phi_{\breve V_0,\breve V_1}([V_0])\in H_n(\breve V_n)$ (here $\breve V_k$ is
the unique $k$-dimensional stratum in $V_k$). In the section \ref{Section Reciprocity Law} we prove
that
$$
res_F \omega=\frac{1}{(2\pi i)^n}\int\limits_{\Delta_F}\omega.
$$
The construction of the Leray coboundary operators is very geometric. In particular, the class
$\Delta_F$ is naturally represented by a smooth submanifold $\tau_F\subset\breve V_n,$ which is a
union of smooth $n$-dimensional tori $\tau_F=\bigcup \tau_{F,a_i}$.

In fact, in the original Parshin's construction, the residue at the flag $F$ is defined as a sum of
certain more delicate residues (we briefly review Parshin's definitions in the Section \ref{Section Parshin's theory reviewed}). We show, that the tori $\tau_{F,a_i}$ naturally correspond to the
summands in the Parshin's definition.

\begin{exam}
Let $S\subset \mathbb C^3$ be the algebraic surface, given by the equation $\{xyz^2+x^4+y^4=0\}.$ Consider the flag $F=\{V_2\supset V_1\supset V_0\},$ where $V_2$ is the surface $S,$\ $V_1$ is the $z$-axis (which is the singular locus of $S$), and $V_0$ is the origin.

The intersection of $S$ with the real plane is the cone over a figure eight (see Figure \ref{figure eight}). The real picture helps a lot to visualize this example.

\begin{figure}[h]
\begin{center}
\includegraphics[width=10cm]{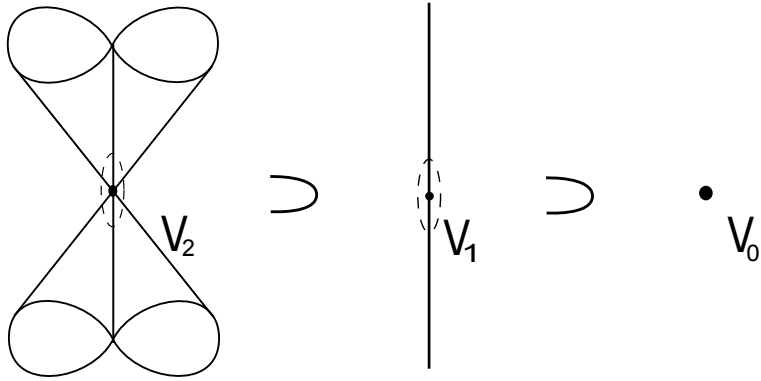}
\caption{This picture shows the intersection of the flag $F$ with the real space.}\label{figure eight}
\end{center}
\end{figure}

There is a natural stratification of $S,$ consisting of $3$ strata: the origin, the $z$-axis without the origin, and the regular part of $S.$ Following the above notations, we denote the strata $\breve V_0,$ $\breve V_1,$ and $\breve V_2$ correspondingly.

$V_0$ is a point on the complex line $V_1.$ One can consider a small circle $\tau^1$ going counterclockwise around $V_0$ on $V_1.$ $\tau^1$ naturally represents the class $\phi_{\breve V_0,\breve V_1}([V_0])\in H_1(\breve V_1).$

On the next step, we have two branches of $V_2$ at each point of $\tau^1.$ Take a point $x\in \tau^1.$ Consider a transversal section to $V_1$ through $x.$ Its intersection with $V_2$ is a curve with two local branches at $x.$ Consider two small circles $S_1$ and $S_2$ around $x,$ one on each branch.

One can choose transversal sections to $V_1$ at each point of $\tau^1,$ so that they depend nicely on the point $x\in \tau^1.$ Furthermore, one can choose the circles in such a way that they form a fiber bundle over $\tau^1.$ Easy to see, that the local branches of $V_2$ do not interchange, as one goes around the origin along $\tau^1$. Therefore, one gets two tori $\tau_{F,a_1}$ and $\tau_{F,a_2}$ in $\breve V_2.$

The Parshin residue of a meromorphic $2$-form $\omega$ on $S$ can be computed via integral over $\tau_F:=\tau_{F,a_1}\cup\tau_{F,a_2}:$
$$
res_F \omega=\frac{1}{(2\pi i)^2}\int\limits_{\tau_F}\omega.
$$
\end{exam}

\begin{exam}
Consider the Whitney umbrella --- the surface $S\subset\mathbb C^3,$ given by the equation $\{y^2-zx^2=0\}.$ Consider the flag $F=\{V_2\supset V_1\supset V_0\},$ where, again, $V_2$ is the surface, $V_1$ is the $z$-axis, and $V_0$ is the origin (Figure \ref{umbrella}). Notice, that $V_1$ again coincide with the singular locus of $S.$

\begin{figure}[h]
\begin{center}
\includegraphics[width=10cm]{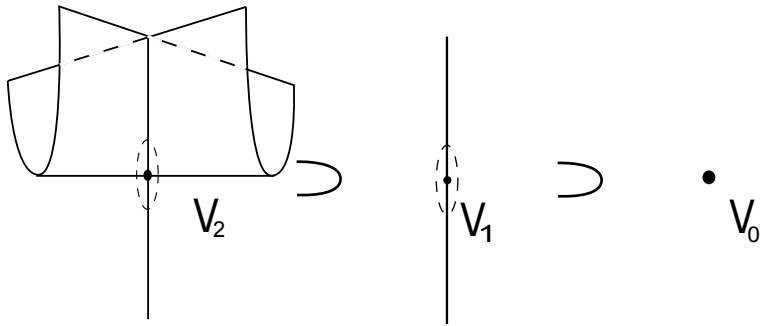}
\caption{As before, this picture shows the intersection of the flag $F$ with the real space.}\label{umbrella}
\end{center}
\end{figure}

In the same way as in the previous example, one can consider a small loop around the origin on the $V_1.$ And, again, $V_2$ has two branches at each point of the loop. However, as one goes around the origin on $V_1,$ the branches interchange. Therefore, the class $F_{\Delta}$ is represented by just one torus in this case, and there is only one summand in the Parshin's definition of the residue.
\end{exam}

Coboundary operators satisfy an interesting relation: Let $X<Y$ be two strata such that there are $k$ intermediate strata $Z_1,\dots,Z_k,$ and these intermediate strata are incomparable (equivalently, for any $m$ from $1$ to $k,$ $X<Z_m<Y$ are consecutive strata). Then
$$
\phi_{Z_1,Y}\circ\phi_{X,Z_1}+\phi_{Z_2,Y}\circ\phi_{X,Z_2}+\dots+\phi_{Z_k,Y}\circ\phi_{X,Z_k}=0.
$$
(see Theorem \ref{Relation}.)

One can illustrate this relation on an example:

\begin{exam}\label{example relation}
On the figure \ref{picture relation}, $X$ is the origin, $Z_1$ is a half-line, $Z_2$ is a surface with an isolated
singularity at the origin, and $Y=\mathbb R^3\backslash(X\cup Z_1\cup Z_2).$ We take a small sphere
$S^2$ with center at the origin. Then $\phi_{X,Z_i}([X])\in H_{\dim Z_i -1}(Z_i)$ is represented by
the intersection $N_i=S^2\cap Z_i.$ Take a small neighborhood of $N_i$ in $S^2.$ Its boundary $D_i$
represents the class $\phi_{Z_i,Y}\circ\phi_{X,Z_i}([X])\in H_1(Y).$ Then the sphere $S^2$ with the
neighborhoods of $N_i$'s deleted gives a two-dimensional chain in $Y,$ which boundary is the union
$D_1\cup D_2.$

\begin{figure}[h]
\begin{center}
\includegraphics[height=7cm]{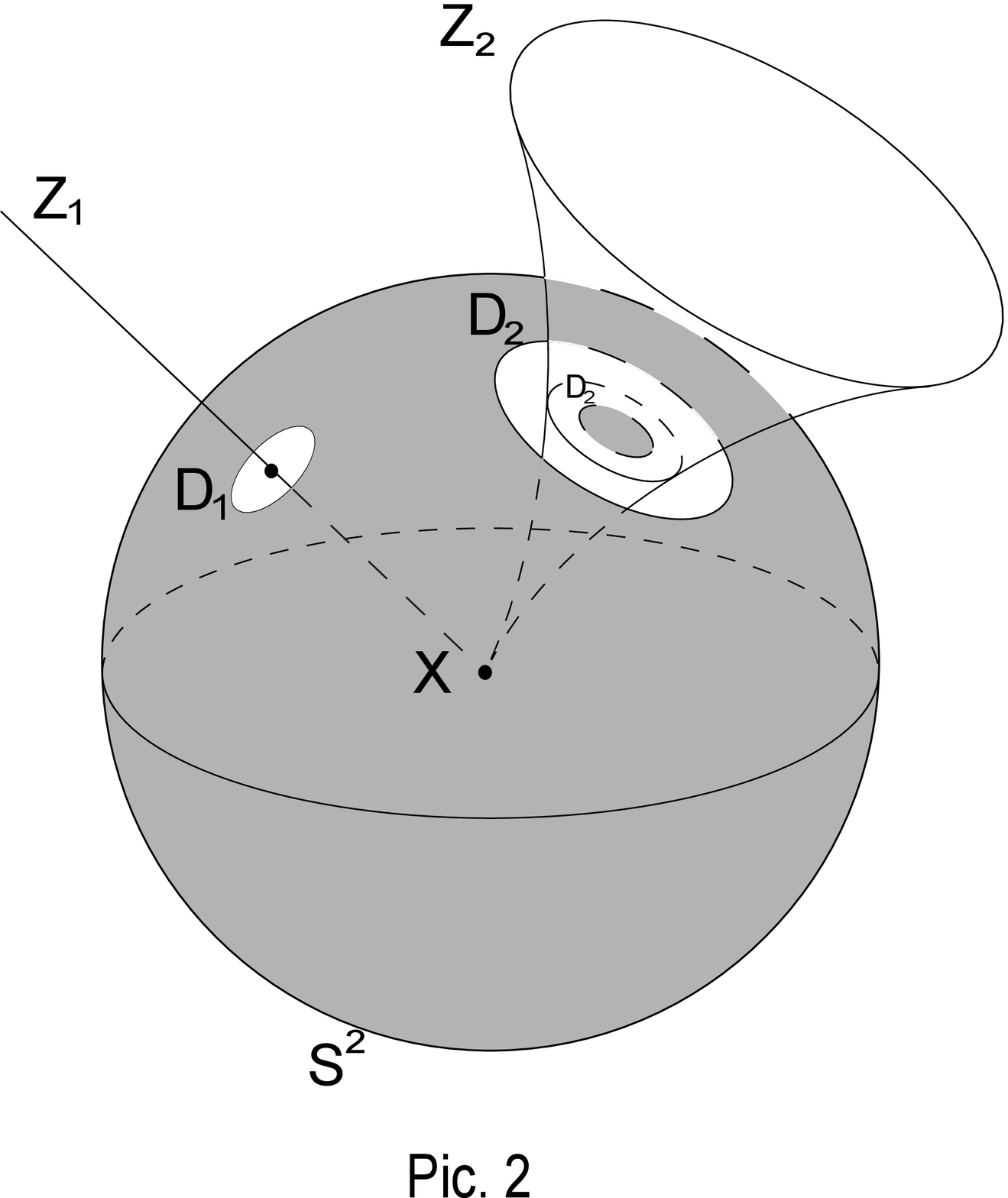}
\caption{Example \ref{example relation}.}\label{picture relation}
\end{center}
\end{figure}
\end{exam}

Our approach also allows us to prove an interesting result about Parshin residues:

Let $\omega$ be a meromorphic top-form on $V_n.$ Consider any Whitney stratification of $V_n$ such
that $\omega$ is regular on the top-dimensional stratum. Then the residue $res_F \omega$ could
be non-trivial, only if all elements of the flag $F$ are closures of strata of the stratification
(Theorem \ref{criteria for non-trivial residues}). In particular, there are only finitely many
non-trivial residues for a given form.

\noindent{\bf Structure of the paper.} In the first part of the paper we introduce the Leray coboundary operators for stratified spaces
and prove the relation (Theorem \ref{Relation}). In the second part we use the results of the first part to express the Parshin residue as an integral over a real smooth cycle and to prove the Reciprocity Law. In the Section \ref{section strat review} we give a short introduction to the theory of stratified spaces. In the Section \ref{Section Parshin's theory reviewed} we review
the original Parshin's definitions and the formulation of the Reciprocity Law.

\noindent {\bf Acknowledgments.} Author would like to thank his PhD advisor, Professor Askold
Khovanskii, for raising the question and helpful discussions.

\section{Leray coboundary operators for stratified spaces.}
\subsection{Whitney Stratifications and Mather's Abstract Stratified Spaces Reviewed.}\label{section strat review}

\begin{df}
Let $M$ be a smooth manifold. Let $V$ be a locally closed subset of $M.$ By a {\it Whitney
stratification} $\bf S$ of $V,$ we mean a subdivision of $V$ into smooth strata, such that:

\begin{enumerate}
\item It is {\it locally finite} - each point of $V$ has an open neighborhood which intersects only
finitely many strata.
\item {\it Condition of the frontier} - for each stratum $X\in\bf S$ its boundary $(\overline X
\backslash X)\cap V$ is a union of strata.
\item  Each pair $(X,Y)$ of strata satisfies {\it Whitney conditions} {\bf a} and {\bf b}:
\begin{enumerate}
\item[{\bf a:}] For any $x \in X$ and any sequence $\{ y_n\} \in Y,$ such that $y_n
\rightarrow x,$ if the sequence of tangent planes $T_{y_n}Y$ converges to some plane $\tau \subset
T_xM$ (in the appropriate Grassmanian bundle over $M$), then $T_xX \subset \tau.$

\item[{\bf b:}] For any $x \in X,$ any sequence $\{ y_n\} \in Y,$ and any sequence $\{ x_n\}
\in X,$ such that $y_n \rightarrow x$ and $x_n \rightarrow x,$ if the sequence of tangent planes
$T_{y_n}Y$ converges to some plane $\tau \subset T_xM,$ and the sequence of secants
$\overline{x_ny_n}$ converges to some line $l$ (in some smooth coordinate system in $M$), then $l
\subset \tau.$
\end{enumerate}
\end{enumerate}
\end{df}

\noindent{\bf Remark.} Actually, condition {\bf b} implies condition {\bf a}, so it is enough to
require condition {\bf b}.

One can prove, that if a pair of strata $(X,Y)$ satisfies condition {\bf b} and $\overline Y \cap X \neq
\emptyset$, then $\dim X < \dim Y.$

\noindent{\bf Notation.} We say, that $X<Y$ if $\overline Y \cap X \neq \emptyset.$ One can see
that this defines a partial order on the set of strata $\bf S.$

\begin{exam}
Consider the surface in $\mathbb C^3$ given by the equation
$y^2+x^3-z^2x^2=0$ (see Figure \ref{whitney example}). The singular locus of the surface coincide with the $z$-axis. Thus, the
$z$-axis and its complement gives a subdivision of the surface in two smooth pieces. Easy to prove
that this pair satisfies the condition {\bf a,} but doesn't satisfy condition {\bf b} at the origin.
Note, that the small neighborhood of the origin looks very different from the neighborhood of any
other point of the $z$-axis.

\begin{figure}[h]
\begin{center}
\includegraphics[height=7cm]{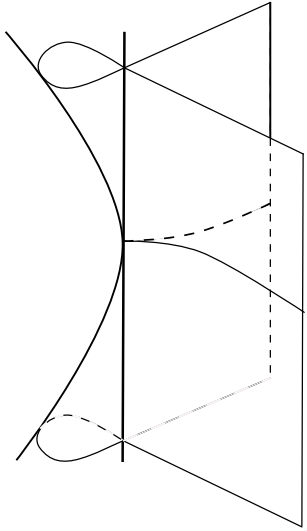}
\caption{As before, the picture shows the intersection of the surface with the real space.}\label{whitney example}
\end{center}
\end{figure}

It is easy to improve the subdivision in such a way that it satisfies condition {\bf b}: one only
needs to consider the origin as a separate stratum.
\end{exam}

Whitney showed that if conditions {\bf a} and {\bf b} are satisfied for the pair $(X,Y)$, then $Y$
"behaves regularly"\ along $X.$

\begin{te}[see \cite{GM}, for example]
Let $V$ be a closed subvariety in a smooth algebraic variety $M.$ Let $\Sigma$ be a locally finite
family of subvarieties in $V.$ Then there exists a Whitney stratification of the set
$V$ such that each element of $\Sigma$ is a union of strata and all strata are algebraic.
\end{te}

Detailed review of the theory of Whitney stratifications can be found in \cite{GM}.

The notion of an Abstract Stratified Space, introduced by John Mather in \cite{MJ}, provides a convenient setup  for working with ``nice'' stratifications: subdivisions into smooth pieces with regular behaviour along strata. J. Mather proved, that any Whitney Stratification can be endowed with a structure of an Abstract Stratified Space. Below we introduce the notion of an Abstract Stratified Space.

Let $V$ be a Hausdorff, locally compact topological space, satisfying the second countability axiom (i.e. there exist a countable basis in the topology of $V$). Let ${\bf S}$ be a locally finite subdivision of $V$ into topological manifolds, endowed with smoothness structures. The elements of $S$ are called strata. Let ${\bf S}$ satisfy the condition of the frontier (i.e. the boundary of any stratum is a union of strata). Similarly as for Whitney Stratifications, the set of strata ${\bf S}$ inherits the natural partial order ($X<Y$ if $X\subset\partial Y$).

For every $X\in{\bf S}$ let $U_X$ be a neighborhood of $X$ in $V,$\  $\rho_X:X\to\mathbb R_{\ge 0}$ be a continuous function, and $\pi_X:U_X\to X$ be a retraction. One should think of $\rho_X$ as of the distance to $X.$ Therefore, we require that $X=\{\rho_X=0\}.$  It is also convenient to say that $\rho_X(y)=\infty$ if $y\notin U_X.$

We call $\rho_X$ --- the tubular function, and $U_X$ --- the tubular neighborhood.

Let $X,Y\in {\bf S},$\ $X\neq Y.$ We use the following notations:

$$
U_{X,Y}:=U_X\cap Y;
$$
$$
\rho_{X,Y}:=\rho_X|_{U_{X,Y}}:U_{X,Y}\to\mathbb R_+;
$$
$$
\pi_{X,Y}:=\pi_X|_{U_{X,Y}}:U_{X,Y}\to X.
$$

We assume that $U_{X,Y}$ is empty unless $X<Y.$ We also assume that if $X$ and $Y$ are incomparable, then $U_X\cap U_Y$ is empty.

We have the following compatibility conditions:

$$
\pi_{X,Y}(\pi_{Y,Z}(v))=\pi_{X,Z}(v),
$$
$$
\rho_{X,Y}(\pi_{Y,Z}(v))=\rho_{X,Z}(v),
$$
\noindent whenever both sides of these equations are defined.

The following conditions ensure that the space $V$ behaves regularly along strata:

For any $X\in {\bf S}$ the map

$$
(\pi_X,\rho_X)|_{U_X\backslash X}:U_X\backslash X \to X\times\mathbb R_+
$$
\noindent is a locally trivial fibration with a compact fiber. Moreover, for any $Y>X$ the restriction

$$
(\pi_{X,Y},\rho_{X,Y}):U_{X,Y}\to X\times\mathbb R_+
$$
\noindent is a smooth fibration.

Finally, we want the fiber of $\pi_X$ over a point $x\in X$ to be a cone with the vertex at $x.$ It does not follow from the above conditions. So, we need to add one more:

Let $U_X^{\le 1}=\{y\in U_X|\rho_X(y)\le 1\}.$ Let $N_X=\partial U_X^{\le 1}=\{y\in U_X|\rho_X(y)=1\}.$ Then $\pi_X|_{U_X^{\le 1}}:U_X^{\le 1}\to X$ is the mapping cone over $\pi_X|_{N_X}:N_X\to X.$

\begin{df}
The Triple ${\bf J}=\{\{U_X\},\{\pi_X\},\{\rho_X\}\}$ is called {\it control data}.
\end{df}

\begin{df}
The Triple $\{V,{\bf S},{\bf J}\}$ with the above conditions is called an {\it Abstract Stratified Space}.
\end{df}

It follows that $N_X$ has a natural structure of an Abstract Stratified Space, obtained by intersecting the strata of $V$ with $N_X$ and restricting tubular functions and retractions.

Abstract Stratified Spaces were introduced by John Mather in \cite{MJ}. In fact, the original definition is slightly different and less restrictive. However, it follows easily that by shrinking the tubular neighborhoods and rescaling the tubular functions, one can change the control data so that they satisfy the above conditions.

\subsection{Leray coboundary operators and relations.}

Let $f:M\to N$ be a smooth fibration with compact oriented $k$-dimensional fiber $F.$ Then one can
define the Gysin homomorphism on homology $f^*:H_*(N)\to H_{*+k}(M).$ Basically, one just set
$f^*(a)=[f^{-1}(A)],$ where $A$ is a representative of the homology class $a\in H_*(N).$

\noindent {\bf Remark.} We use the following convention about the orientations. Let $y\in M$ and
$x=f(y)\in N.$ Let $A\subset N$ be a smooth representative of a homology class $a\in H_*(N)$ and
$x\in A.$ Let the differential form $\omega_A$ on $N$ be such that its restriction to $A$ defines
the orientation of $A$ at $x$ and the differential form $\omega_F$ on $M$ be such that its
restriction to the fiber $F_x$ defines the orientation of $F_x$ at $y.$ Then the orientation of the
preimage $f^{-1}(A)\subset M$ at the point $y$ is given by the restriction of the form
$f^*(\omega_A)\wedge\omega_F.$

Let now $M$ be an oriented manifold with boundary, and $f:M\to N$ be a proper map to an oriented
manifold $N$, such that its restriction both to the boundary $\partial M\subset M$ and the interior
$\breve M\subset M$ are submersions. Then, by the Ehresmann Lemma for manifolds with boundary, $f$
is a locally trivial fibration and its restrictions to $\partial M$ and $\breve M$ are smooth
fibrations.

Let $\phi:=(f|_{\partial M})^*:H_*(N)\to H_{*+\dim M-\dim N-1}(\partial M)$ be the Gysin
homomorphism.

\begin{lem}\label{Relation for Gysin maps}
$i_*\circ\phi=0,$ where $i:\partial M\hookrightarrow M$ is the embedding.
\end{lem}

\begin{proof}
One can generalize the Gysin homomorphism to the described above case, when the fiber of $f$ is a
manifold with boundary. The only difference is that now the homomorphism lands in the relative
homology group: $f^*:H_*(N)\to H_{*+k}(M,\partial M),$ where $k=\dim F=\dim M-\dim N.$ Then one
immediately sees that $\phi=\partial\circ f^*,$ where $\partial:H_*(M,\partial M)\to
H_{*-1}(\partial M)$ is the boundary homomorphism from the long exact sequence of the pair
$(M,\partial M).$ However, by the long exact sequence, $i_*\circ\partial=0.$
\end{proof}

We apply the above constructions to the stratified spaces.

Let all the strata of a stratified space $V$ be oriented.

Let $X\in\bf S$ be a stratum. Restriction of the retraction
$\pi_X:U_X\to X$ to $N_X=\{y\in U_X| \rho_X(y)=1\}$ is a locally trivial fibration. Moreover, for
any stratum $Y$ such that $X<Y$ the restriction to $N_{X,Y}=N_X\cap Y=\{y\in U_{X,Y}|
\rho_X(y)=1\}$ is a smooth fibration.

\begin{df}\label{order on strata}
Let $X<Y$ be two strata. We say that $X<Y$ are {\it consecutive strata} if there is no such $Z$
that $X<Z<Y.$
\end{df}

\begin{lem}
Let $X<Y$ be consecutive strata. Then the fiber of $\pi_X|_{N_{X,Y}}:N_{X,Y}\to X$ is compact.
\end{lem}

\begin{proof}
Since $X<Y$ are consecutive strata, it follows that $N_{X,Y}=Y\cap N_X$ is a closed stratum of
$N_X$ (indeed, otherwise the closure of $N_{X,Y}$ in $N_X$ would contain a smaller stratum). The
fiber of the restriction of $\pi_X$ to $N_{X,Y}$ is the intersection of the fiber of the
restriction of $\pi_X$ to $N_X$ and $N_{X,Y}.$ Therefore, it is compact as a closed subset of a
compact set.
\end{proof}

Note, that $N_{X,Y}$ is orientable. Indeed, it is the level set of a smooth function $\rho_{X,Y}$
in $U_{X,Y}\subset Y.$ Let us fix the orientation of $N_{X,Y}$ given as follows: we say that the
restriction a differential $(\dim Y -1)$-form $\omega_{N_{X,Y}}$ on $Y$ defines the positive
orientation of $N_{X,Y}$ if the form $d\rho_{X,Y}\wedge \omega_{N_{X,Y}}$ defines the positive
orientation of $Y.$

Let $\dim X=n$ and $\dim Y=k.$

\begin{df}
The {\it Leray coboundary operator} $\phi_{X,Y}:H_*(X)\to H_{*+k-n-1}(Y)$ is given by the
composition $\phi_{X,Y}=i_*\circ\phi,$ where $i:N_{X,Y}\hookrightarrow Y$ is the embedding and
$\phi:H_*(X)\to H_{*+k-n-1}(N_{X,Y})$ is the Gysin homomorphism.
\end{df}

\begin{te}\label{Relation}
Let $X<Y$ be two strata. Let $Z_1,\dots,Z_m$ be all strata such that $X<Z_i<Y.$ Suppose that
$Z_1,\dots,Z_m$ are incomparable. Then
$$
\phi_{Z_1,Y}\circ\phi_{X,Z_1}+\phi_{Z_2,Y}\circ\phi_{X,Z_2}+\dots+\phi_{Z_m,Y}\circ\phi_{X,Z_m}=0.
$$
\end{te}

\begin{proof}
We want to apply the Lemma \ref{Relation for Gysin maps}. Consider $D_i:=N_{X,Y}\cap
N_{Z_i,Y}=\{y\in Y|\rho_{Z_i}(y)=\rho_X(y)=1\}.$ Note, that
$D_i=(\pi_{Z_i}|_{N_{Z_i,Y}})^{-1}(N_{X,Z_i}).$ Therefore, $\pi_{Z_i}|_{D_i}$ is a smooth fibration
over $N_{X,Z_i}.$ Denote $p_i:=\pi_X|_{N_{X,Z_i}}\circ\pi_{Z_i}|_{D_i}: D_i\to X.$ According to the
construction of the Leray coboundary operators, we have
$$
\phi_{Z_i,Y}\circ\phi_{X,Z_i}=i_*\circ \phi_i,
$$
\noindent where $i:D_i\hookrightarrow Y$ is the embedding and $\phi_i:H_*(X)\to H_{*+\dim Y-\dim
X-2}$ is the Gysin homomorphism of $p_i:D_i\to X.$ Here we fix the orientation of $D_i$ given in
the following way: we say that the restriction of a differential $(\dim Y -2)$-form $\omega_{D_i}$
on $Y$ defines the positive orientation of $D_i$ if the form $d\rho_{Z_i,Y}\wedge d\rho_{X,Y}\wedge
\omega_{D_i}$ defines the positive orientation of $Y.$

Consider now $N_{X,Y}=\{y\in Y|\rho_X(y)=1\}.$ The restriction $\pi_X|_{N_{X,Y}}$ is a smooth
fibration. However, the fibers of this fibration are not compact. On the other side, if we consider
the restriction of $\pi_X$ to the union $N_{X,Y\cup Z_1\cup\dots\cup Z_m}:=N_{X,Y}\cup
N_{X,Z_1}\cup\dots\cup N_{X,Z_m}= N_X\cap (Y\cup Z_1\cup\dots\cup Z_m),$ then the fibers are
compact.

$D_i\subset N_{X,Y}$ can be thought of as the boundary of the neighborhood $U_i=\{y\in N_{X,Y}\cap
U_{Z_i}| \rho_{Z_i}(y)<1\}$ of $N_{X,Z_i}$ in $N_{X,Y\cup Z_1\cup\dots\cup Z_m}.$ Denote
$M=N_{X,Y}\backslash (U_1\cup\dots\cup U_m).$ By Ehresmann Lemma for manifolds with boundary, the
restriction $\pi_X|_M:M\to X$ is a locally trivial fibration. Indeed, $\pi_X|_M$ is proper, because
$M$ is a closed subset of $N_{X,Y\cup Z_1\cup\dots\cup Z_m}$ and $\pi_X|_{N_{X,Y\cup
Z_1\cup\dots\cup Z_m}}$ is a fibration with compact fibers; the restrictions of $\pi_X$ to the
interior of $M$ and the boundary $\partial M=D_1\cup\dots\cup D_m$ are submersions.

To conclude the proof by Lemma \ref{Relation for Gysin maps}, one needs to check, that the
orientation of $D_i$ as a piece of the boundary of $M$ always coincide (or always is opposite) with the
orientation of $D_i$ used in the first part of the proof. Indeed, we fixed the orientation of $D_i$
in such a way, that if $\omega_{D_i}|_{D_i}$ gives the orientation of $D_i$ then
$d\rho_{Z_i,Y}\wedge d\rho_{X,Y}\wedge \omega_{D_i}$ gives the orientation of $Y.$ Let
$\omega_{N_{X,Y}}:=-d\rho_{Z_i,Y}\wedge \omega_{D_i}.$ According to our convention about the
orientation of $N_{X,Y},\ $ $\omega_{N_{X,Y}}$ gives the positive orientation of $N_{X,Y}.$
Therefore, the orientation of $D_i$ as a piece of the boundary of $M$ is given by $-\omega_{D_i}.$
\end{proof}

\subsection{Dual Homomorphism.}

In this chapter the coefficient ring is always $\mathbb R$. For simplicity, we skip it in the
notations.

\smallskip

There is a natural question:

{\it Which operator is Poincare dual to the coboundary operator $\phi_{X,Y}?$}

The manifolds $X$ and $Y$ are not compact. Therefore, one has to use the Borel-Moore homology to do the Poincare duality. A nice review of the theory of Borel-Moore homology (in much more details than needed here) is given in \cite{Gi}.

Consider $\phi_{X,Y}: H_m(X)\to H_{m+k-n-1}(Y)$ (here $\dim X=n,$ $\dim Y=k$). The dual operator is $(\phi_{X,Y})^*:H^{BM}_{n-m+1}(Y)\to H^{BM}_{n-m}(X).$ There is a natural candidate for the dual: indeed, one can show, that

$$
H_{n-m+1}^{BM}(Y)=H_{n-m+1}^{BM}(Y\cup X,X).
$$

Therefore, there exits the boundary operator

$$
\partial_{Y,X}:H^{BM}_{n-m+1}(Y)\to H^{BM}_{n-m}(X).
$$

{\bf Remark.} It is crucial that $X<Y$ are consecutive strata. Otherwise, the union $X\cup Y$ would not be locally compact, and the boundary operator would not be defined.

\begin{te}
Leray coboundary operator $\phi_{X,Y}:H_m(X)\to H_{m+k-n-1}(Y)$ ($\dim X=n,$ $\dim Y=k$) is
Poincare dual to the boundary homomorphism $\partial_{Y,X}:H^{BM}_{n-m+1}(Y)\to H^{BM}_{n-m}(X).$
\end{te}

\begin{proof}
By Poincare duality, the intersection form $H_*(M)\times H^{BM}_{d-*}(M)\to \mathbb R$ is well
defined and non-degenerate (here $M$ is a smooth oriented manifold and $\dim M=d$). Therefore, the
only thing we need to check is that for any classes $a\in H_n(X)$ and $b\in H_{m-n+1}^{BM}(Y),$

$$
<\partial_{Y,X} b,a>=<b,\phi_{X,Y}(a)>.
$$

Let $i:U_{X,Y}\hookrightarrow Y$ be the embedding. According to the definition of the Leray
coboundary operator, $\phi_{X,Y}$ can be factored: $\phi_{X,Y}=i_*\circ\phi_{X,U_{X,Y}},$ where
$\phi_{X,U_{X,Y}}:H_m(X)\to H_{m+k-n-1}(U_{X,Y})$ is the Leray coboundary operator for the
stratified space with two strata: $X$ and $U_{X,Y}.$ On the other side, the boundary homomorphism
$\partial_{Y,X}$ also can be factored: $\partial_{Y,X}=\partial_{U_{X,Y},X}\circ i^*$ (here
$i^*:H^{BM}(Y)\to H^{BM}(U_{X,Y})$ is the restriction homomorphism induced by the inclusion $i$). Therefore, it
is enough to assume that $Y=U_{X,Y}.$

We know that $U_{X,Y}$ is diffeomorphic to
$N_{X,Y}\times\mathbb R_+.$ Therefore, there is an isomorphism
$\theta:H_*^{BM}(U_{X,Y})\stackrel{\sim}{\longrightarrow} H_{*-1}^{BM}(N_{X,Y}),$ given by taking a
representative, transversal to $N_{X,Y}$ and intersecting it with $N_{X,Y}.$ The inverse
isomorphism $\theta^{-1}$ is given by multiplying a representative by $\mathbb R_+.$

{\bf Remark.} One should be careful with the orientations. We want the following
condition to be satisfied: if $A\subset U_{X,Y}$ is a cycle, transversal to $N_{X,Y},\ $
$B:=N_{X,Y}\cap A,$ and $\omega_B|_B$ gives the orientation of $B$ at some point, then
$d\rho_{X,Y}\wedge\omega_B$ should give the positive orientation of $A$ at this point.

With the above orientation conventions one gets that
$$
<b,\phi_{X,Y}(a)>=<\theta(b),\phi(a)>,
$$
\noindent where $\phi:H_*(X)\to H_{*+k-n-1}(N_{X,Y})$ is the Gysin homomorphism and the
intersection on the right is taken inside $N_{X,Y}$. Moreover,

$$
\partial_{U_{X,Y},X}=(\pi_X|_{N_{X,Y}})_*\circ\theta.
$$

Therefore, we only need to check, that the Gysin homomorphism is dual to the
$(\pi_X|_{N_{X,Y}})_*:H_{n-m}(N_{X,Y})\to H_{n-m}(X),$ which is obvious.

\end{proof}

\noindent{\bf Corollary.} Leray coboundary operator $\phi_{X,Y}$ does not depend on the choice of
the control data at least modulo torsion.

One can also investigate the relation, dual to the relation on the coboundary operators, proved in the Theorem \ref{Relation}.

Consider the strata $X,Z_1,\dots,Z_p,Y,$ satisfying the conditions of the Theorem \ref{Relation}. Let $Z=\bigcup Z_i.$ The boundary operator

$$
\partial_{Y,Z}:H_*^{BM}(Y)\to H_{*-1}^{BM}(Z)=\oplus H_{*-1}^{BM}(Z_i)
$$

\noindent is dual to the direct sum of the coboundary operators $\oplus \phi_{Z_i,Y}.$ The boundary operator

$$
\partial_{Z,X}:H_{*-1}^{BM}(Z)\to H_{*-2}^{BM}(X)
$$

\noindent is, in turn, dual to the direct sum $\oplus \phi_{X,Z_i}.$ Therefore, the dual relation is

$$
\partial_{Z,X}\circ\partial_{Y,Z}=0.
$$

However, it is not hard to prove independently. Indeed, one can use that $H_*^{BM}(Y)=H_*^{BM}(Y\cup Z\cup X,Z\cup X)$ and $H_*^{BM}(Z)=H_*^{BM}(Z\cup X,X).$ Then the relation basically says, that the boundary of the boundary of a chain is zero, which is trivial. This provides another proof of the Theorem \ref{Relation} modulo torsion.

{\bf Remark.} It is crucial for this argument, that $Y\cup Z\cup X$ is locally compact. Therefore, if one would accidentally forget one of the intermediate strata $Z_i,$ the relation would not hold.

%One can consider the spectral sequence for the
%Borel-Moore homologies of this filtration. Note, that all three terms of the filtration are locally
%compact. Therefore, one can easily show that the first term of the spectral sequence is given by

%$ E_{0,i}^1=H_i^{BM}(X), $

%$ E_{1,i}^1=H_{i+1}^{BM}(X\cup Z_1\cup\dots\cup Z_p,X)=H_{i+1}^{BM}(Z_1)\oplus\dots\oplus
%H_{i+1}^{BM}(Z_p)), $

%$ E_{2,i}^1=H_{i+2}^{BM}(X\cup Z_1\cup\dots\cup Z_p\cup Y,X\cup Z_1\cup\dots\cup
%Z_p)=H_{i+2}^{BM}(Y), $

%$
%\partial_{1,*}^1=\bigoplus_j\partial_{Z_j,X},
%$

%$
%\partial_{2,*}^1=\bigoplus_j\partial_{Y,Z_j}.
%$

%$$
%\begin{array}{|l|l|l|l|l|l|}
%\hline
%H_0^{BM}(Y) & H_1^{BM}(Y) & H_2^{BM}(Y) & H_3^{BM}(Y) & H_4^{BM}(Y) & \dots \\
%\hline
%0 & \bigoplus H_0^{BM}(Z_i) & \bigoplus H_1^{BM}(Z_i) & \bigoplus H_2^{BM}(Z_i) & \bigoplus H_3^{BM}(Z_i) & \dots \\
%\hline
%0 & 0 & H_0^{BM}(X) & H_1^{BM}(X) & H_2^{BM}(X) & \dots \\
%\hline
%\end{array}
%$$

%Now, the condition that the square of the boundary operator is zero is dual to the relation on coboundary operators, given by the Theorem \ref{Relation}.

%This provides another proof to the Theorem \ref{Relation} modulo torsion.

\section{Application to Parshin's Residues.}
\subsection{Parshin Residues and the Reciprocity Law Reviewed.}\label{Section Parshin's theory reviewed}

In this section we review the definition of the Parshin residue and the Reciprocity Law. As it was discussed in the introduction, Parshin residue at a flag $F$ is defined as a sum of certain more delicate residues. In fact, every flag ``contains'' finitely many {\it Parshin points}, and the more delicate residues are computed at these points.

We start from the definition of a Parshin point:

Let $V_n$ be an algebraic variety of dimension $n.$ Let $F=:\{V_n\supset\dots\supset V_0\}$ be a flag of subvarieties of dimensions $\dim V_k=k.$

Consider the following diagram:

\begin{equation}\label{NormDiag}
\begin{CD}
V_n             @.\supset @. V_{n-1}                     @.\supset @.  \dots @.\supset @. V_1                            @.\supset @. V_0\\
@AAp_nA              @.      @AAp_nA                          @.        @.       @.             @.                           @.       @. \\
\widetilde{V}_n @.\supset @. W_{n-1}                     @.        @.       @.         @.                                @.        @.    \\
@.                   @.      @AAp_{n-1}A                      @.        @.       @.             @.                           @.       @. \\
                @.        @. \widetilde{W}_{n-1}         @.\supset @. \dots @.        @.                                @.        @.   \\
@.                   @.      @.                               @.        @.       @.             @.                           @.       @. \\
                @.        @.                             @.        @. \dots @.\supset @. W_1                            @.        @.   \\
@.                   @.      @.                               @.        @.       @.       @AAp_{1}A                          @.       @. \\
                @.        @.                             @.        @.       @.        @. \widetilde{W}_1                @.\supset @.  W_0
\end{CD}
\end{equation}

\noindent where
\begin{enumerate}
\item $p_n:\widetilde{V}_n\to V_n$ is the normalization;
\item $W_{n-1}\subset\widetilde{V}_n$ is the union of $(n-1)$-dimensional irreducible components of the preimage of $V_{n-1};$
\item for every $k=1,2\dots,n-1$
\begin{enumerate}
\item $p_k:\widetilde{W_k}\to W_k$ is the normalization;
\item $W_{k-1}\subset\widetilde{W}_k$ is the union of $(k-1)$-dimensional irreducible components of the preimage of $V_{k-1}.$
\end{enumerate}
\end{enumerate}

\begin{df}
We call diagram \ref{NormDiag} the {\it normalization diagram} of the flag $V_n\supset\dots\supset
V_0.$
\end{df}

\begin{df}\label{Definition Parshin point}
The flag $F=\{V_n\supset\dots\supset V_0\}$ of irreducible subvarieties together with a choice of a
point $a_{\alpha}\in W_0$ is called a {\it Parshin point}.
\end{df}

Choosing a point $a_{\alpha}\in W_0$ is equivalent to choosing irreducible components in every
$W_i,$ $i=n-1,\dots, 0.$ Indeed, $\widetilde{W}_i$ is normal and, therefore, locally irreducible at
every point. In particular, it is locally irreducible at the image of $a_{\alpha}.$ Let
$\widetilde{W}_i^{\alpha}$ be the irreducible component of $\widetilde{W}_i,$ containing the image
of $a_{\alpha}.$ Let $W_i^{\alpha}=p_i(\widetilde{W}_i^{\alpha}).$ Note, that $W_i^{\alpha}$ is an
irreducible component of $W_i.$

\begin{exam}
Consider the flag from the Example \ref{figure eight} in the Introduction. The normalization diagram looks as follows:
%\newpage
\begin{figure}[h]
\begin{center}
\includegraphics[width=6cm]{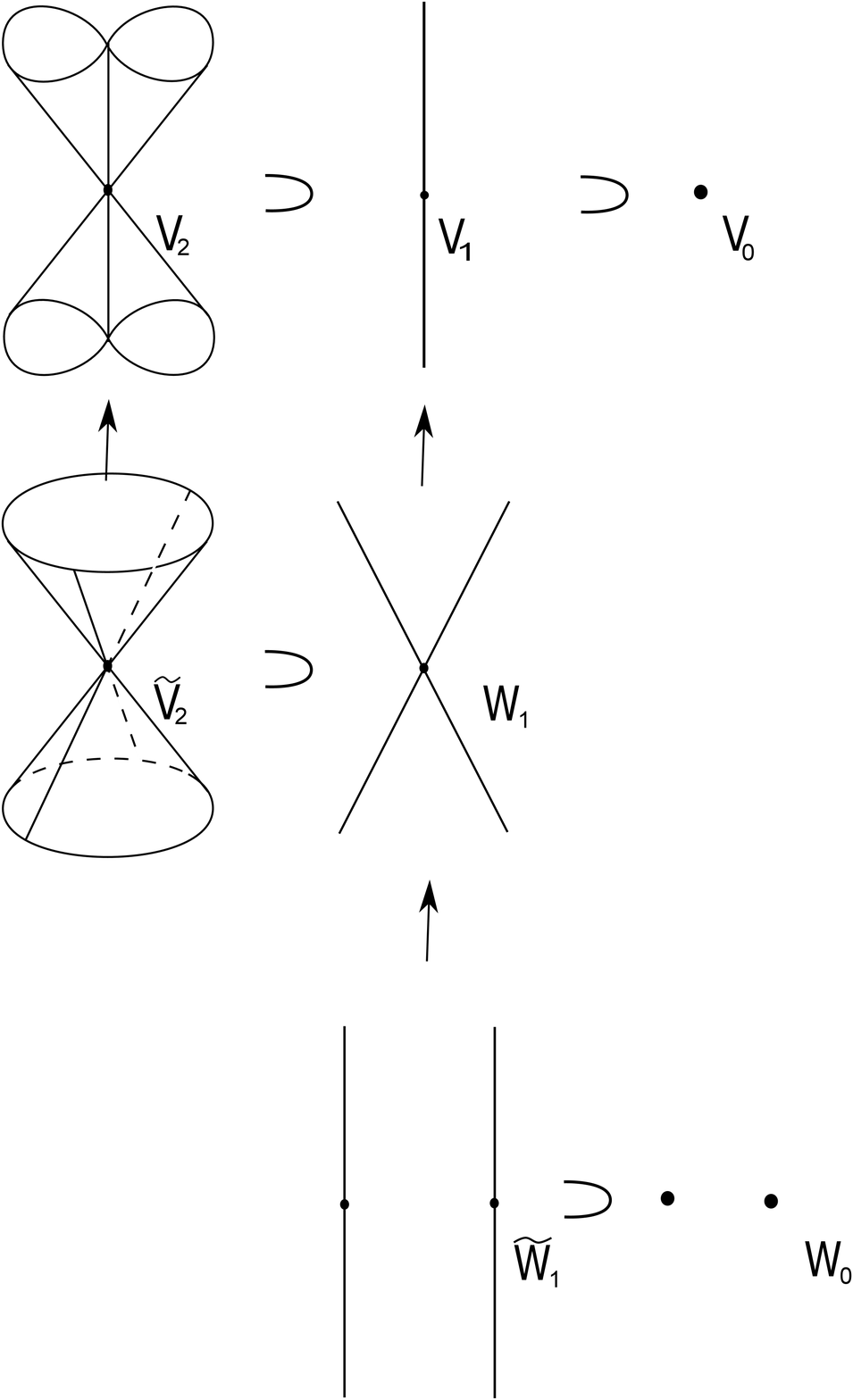}
\end{center}
\caption{\small As in the introduction, we draw the intersection with the real space. Normalization splits the local irreducible components at every point. Therefore, the normalization of the cone over the figure eight is the usual cone, and the preimage of the $z$-axis is two lines, intersecting at the origine.} \label{Norm Diagram figure eight}
\end{figure}

Therefore, there are two Parshin points, corresponding to the flag. Note, that these points naturally correspond to the tori from the Example \ref{figure eight}.
\end{exam}

\begin{exam}
The normalization of the Whitney umbrella (Example \ref{umbrella} in the Introduction) is isomorphic to $\mathbb C^2.$ The preimage of the $z$-axis is a line, which covers the $z$-axis twice with a branching at the origin. Therefore, $W_0$ is just one point. This corresponds to the fact that there is only one torus in the Example \ref{umbrella}.
\end{exam}

In order to define the Parshin residue, one needs to define the {\it local parameters} at a Parshin
point, which play the role of the normalizing parameter in one-dimensional case. After that, one
uses these parameters to define a sequence of residual meromorphic forms
$\omega_{n-1},\dots,\omega_0$ on $W^{\alpha}_{n-1},\dots,W^{\alpha}_0$.

The local parameters are defined as follows:

$W_{i-1}^{\alpha}\subset \widetilde{W}_i$ is a hypersurface in a normal variety. It follows that
there exists a (meromorphic) function $u_i$ on $\widetilde{W}_i$ which has zero of order $1$ at a
generic point of $W_{i-1}^{\alpha}.$ Since meromorphic functions are the same on $W_i$ and
$\widetilde{W}_i,$ one can consider $u_i$ as a function on $W_i.$ Then one can extend (in an
arbitrary way) $u_i$ to $\widetilde{W}_{i+1}$ and so on. For simplicity, we denote all these
functions by $u_i.$ Now $u_i$ is defined on $V_n,$ and can be consecutively restricted to $W_j$ for
$j\ge i.$

\begin{df}\label{Definition parameters}
Functions $(u_1,\dots,u_n)$ are called {\it local parameters} at the Parshin point
$P=\{V_n\supset\dots\supset V_0,a_{\alpha}\}$.
\end{df}

\noindent {\bf Remark.} One can choose local parameters in such a way, that $u_i$ has zero of order
$1$ at generic point not only of $W_{i-1}^{\alpha},$ but of the whole $W_{i-1}.$ Then these
local parameters work for all Parshin points with the flag $F=\{V_n\supset\dots\supset V_0\}.$ We
only use local parameters with this property.

Let $\omega$ be a meromorphic $n$-form on $V_n.$ One can show that the differentials
$du_1,\dots,du_n$ are linearly independent at a generic point of $V_n.$ Therefore, one can write

$$
\omega=fdu_1\wedge\dots\wedge du_n,
$$

\noindent where $f$ is a meromorphic function on $V_n.$

Now we define the residual forms $\omega_i$:

Take a generic point $p\in W^{\alpha}_{n-1}.$ Both $\widetilde{V}_n$ and $W_{n-1}$ are smooth at
$p.$ Moreover, parameters $u_1,\dots, u_n$ provide an isomorphism of a neighborhood of $p\ $ to an
open subset in $\mathbb C^n,$ and $W^{\alpha}_{n-1}$ is given by the equation $u_n=0$ in this
neighborhood. Restrict the function $f$ to the transversal section to $W_{n-1}$ at $p,$ given by
fixing the parameters $u_1,\dots, u_{n-1}.$ The restriction can be expanded into a Laurent series
in $u_n.$ It is easy to see that the coefficients of this expansion depend analytically on $p.$
Moreover, one can see that the coefficients are meromorphic functions on $W^{\alpha}_{n-1}.$ Let
$f_{-1}$ be the coefficient at $u_n^{-1}$ in this expansion. Then
$\omega_{n-1}=f_{-1}du_1\wedge\dots\wedge du_{n-1}$ is a meromorphic $(n-1)$-form on
$W^{\alpha}_{n-1}.$

Repeating this procedure one more time one gets a meromorphic $(n-2)$-form on $W^{\alpha}_{n-2}.$
Finally, after $n$ steps, one gets a function $\omega_0$ on the one-point set
$W^{\alpha}_0=\{a_{\alpha}\}.$

\begin{df}
The {\it residue of $\omega$ at the Parshin point $P=\{V_n\supset\dots\supset V_0,a_{\alpha}\in
W_0\}$} is $res_P(\omega)=\omega_0(a_{\alpha}).$
\end{df}

Parshin proves that the residue is independent on the choice of local parameters.

\begin{df}
The sum of residues over all $a\in W_0$ is called the {\it residue at the flag}
$F=\{V_n\supset\dots\supset V_0\}$ and is denoted $res_F(\omega)=\sum\limits_{a\in W_0}
res_{\{F,a\}}(\omega).$
\end{df}

\begin{te}(\cite{P1},\cite{B},\cite{L})\label{Parshin's Reciprocity Law}
Let $\omega$ be a meromorphic $n$-form on $V_n.$ Fix a partial flag of irreducible subvarieties
$\{V_n\supset\dots\supset\widehat{V_k}\supset\dots\supset V_0\}$, where $V_k$ is omitted ($0<k<n$).
Then
$$
\sum\limits_{V_{k+1}\supset X\supset V_{k-1}} res_{V_n\supset\dots\supset X\supset\dots\supset
V_0}(\omega)=0,
$$
\noindent where the sum is taken over all irreducible $k$-dimensional subvarieties $X,$ such that
$V_{k-1}\supset X\supset V_{k+1}.$ (In this formula only finitely many summands are non zero.)

In addition, if $V_1$ is compact then one has the same relation for $k=0.$
\end{te}

\subsection{Residues via Leray Coboundary Operators and the Reciprocity Law.}\label{Section Reciprocity
Law}

We want to apply the stratification theory to study the Parshin points and residues. Therefore, we
need to stratify all the spaces in the normalization diagram in such a way that the stratifications
respect the normalization maps $p_1,\dots,p_n.$ The following Lemma easily follows from the well
known results on existence of Whitney stratifications (see Section 1.7 in \cite{GM}, for example):

\noindent{\bf Notation.} Let $X$ be an irreducible (complex analytic) variety considered with a
fixed Whitney stratification. Then by $\breve X$ we denote the stratum of maximal dimension. If $X$
is reducible, then by $\breve X$ we denote the union of strata of maximal dimension.

\begin{lem}\label{Stratification of the normalization diagram}
Fix a Parshin point $P=\{V_n\supset\dots\supset V_0,a_{\alpha}\in W_0\}$ and local parameters
$u_1,\dots, u_n.$ There exist Whitney stratifications ${\bf S},{\bf S}_{\widetilde{V}},{\bf
S}_{\widetilde{W}_{n-1}},\dots,{\bf S}_{\widetilde{W}_1}$ of
$V_n,\widetilde{V}_n,\widetilde{W}_{n-1},\dots,\widetilde{W}_1$ correspondingly, such that:
\begin{enumerate}
\item $V_{n-1},\dots,V_0$ are unions of strata of ${\bf S};$
\item $W_{n-1},W_{n-2},\dots, W_0$ are unions of strata of ${\bf S}_{\widetilde{V}},{\bf
S}_{\widetilde{W}_{n-1}},\dots,{\bf S}_{\widetilde{W}_1}$ correspondingly;
\item for all $i=1,\dots, n,$ the local parameter $u_i$ is regular and non-vanishing on $\breve V_n,\breve{\widetilde{V}}_n,\breve{\widetilde{W}}_{n-1},\dots,\breve{\widetilde{
W}}_i;$
\item for all $i=1,\dots, n,$ the restriction of the normalization map $p_i$ to any stratum in the source is a covering over a stratum in the image.
\end{enumerate}
\end{lem}

There is an important corollary about stratifications ${\bf S_{\widetilde{W}_i}}:$

\begin{lem}\label{(n-1)-dimensional stratum}
The stratum (or the union of strata if $W_{i-1}$ is reducible) $\breve{W}_{i-1}\in {\bf S_{\widetilde{W}_i}}$ consists of regular points of
$\widetilde{W}_i.$
\end{lem}

\begin{proof}
Let $x \in \breve{W}_{i-1}$ be a point, such that $\widetilde{W}_i$ is singular at $x.$ Note, that
by dimension reasons and condition of the frontier, the only strata intersecting a small
neighborhood of $x$ are  $\breve{W}_{i-1}$ and $\breve{\widetilde{W}_i}.$ Note, that $u_i$
is regular in $\breve{\widetilde{W}}_i$ and at a generic point of $\breve{W}_{i-1}.$ Therefore, by
extension theorem for normal varieties, $u_i$ is regular at $x.$

Note also, that $u_i$ is non-vanishing in $\breve{\widetilde{W}_i}$ and has zero of order $1$ at a
generic point of $W_{i-1}.$ Therefore, $\{u_i=0\}$ coincide with $W_{i-1}$ near $x.$ Moreover, easy
to see, that the germ of $u_i$ at $x$ generates the ideal of the germ of $W_{i-1}$ at $x.$ Indeed,
if $g$ is a function, regular at $x$ and vanishing on $W_{i-1},$ then $\frac{g}{u_i}$ is regular at
$x$ by the extension theorem for normal varieties.

Now, let $f_1,\dots,f_{i-1}$ be any coordinate system on $W_{i-1}$ at $x.$ Easy to see, that the
functions $u_i,f_1,\dots,f_{i-1}$ generate the maximal ideal in the local ring of
$\{x\}\subset\widetilde{W}_i.$ Therefore, $x$ is a smooth point of $\widetilde{W}_i.$
\end{proof}

Our goal is to show, that
$$
res_F(\omega)=\frac{1}{(2\pi i)^n}\int\limits_{\Delta_F}\omega,
$$
\noindent where $F:=\{V_n\supset\dots\supset V_0\}$ and $\Delta_F=\phi_{\breve V_{n-1},\breve
V_n}\circ\dots\circ\phi_{\breve V_0,\breve V_1}([V_0])\in H_n(\breve V_n).$

Moreover, we will show that $\Delta_F$ naturally splits into the sum $\Delta_F=\sum\limits_{a_i\in
W_0} \Delta_{\{F,a_i\}},$ such that
$$
res_{\{F,a_i\}}(\omega)=\frac{1}{(2\pi i)^n}\int\limits_{\Delta_{\{F,a_i\}}}\omega.
$$
Note, that according to the construction of the Leray coboundary operator, $\Delta_F$ is
represented by a smooth compact real $n$-dimensional submanifold $\tau_F\subset\breve V_n.$
Moreover, $\tau_F$ is obtained from a point by the following procedure: there are $n$ steps, and on
each step we take the total space of an oriented fibration with $1$-dimensional compact fiber over
the result of the previous step. Therefore, $\tau_F$ is a union of $n$-dimensional tori. We'll
show, that the connected components of $\tau_F$ are in natural one-to-one correspondence with the
points of $W_0$ and the connected component $\tau_{F,a_i}$ corresponding to $a_i\in W_0$ represents
$\Delta_{F,a_i}.$

Fix control data on the stratification ${\bf S}$ of $V_n.$ Let us use these control data to
construct the representative $\tau_F\subset \breve V_n$ of $\Delta_F.$ Let us also denote by
$\tau_k\subset \breve V_k$ the representative of $\Delta_k=\phi_{\breve V_{k-1},\breve
V_k}\circ\dots\circ\phi_{\breve V_0,\breve V_1}([V_0])\in H_k(\breve V_k),$ constructed in the same
way.

Let us introduce the following notations:

\begin{enumerate}
\item $\widehat{U}_0:={\breve V_0};$
\item $\widehat{U}_k:=\pi_{{\breve V_{k-1}},{\breve V_k}}^{-1}(\widehat{U}_{k-1}),$ for $k=1,\dots,n.$
\end{enumerate}

Note, that for $k>0,\ $ $\widehat{U}_k$ is the preimage of $\widehat{U}_{k-1}\times\mathbb R_+$
under the map $(\pi_{{\breve V_{k-1}},{\breve V_k}},\rho_{{\breve V_{k-1}},{\breve
V_k}}):U_{{\breve V_{k-1}},{\breve V_k}}\to {\breve V_{k-1}}\times\mathbb R_+.$ Since ${\breve
V_{k-1}}$ and ${\breve V_k}$ are consecutive strata, it follows that the
restriction $(\pi_{{\breve V_{k-1}},{\breve V_k}},\rho_{{\breve V_{k-1}},{\breve
V_k}})|_{\widehat{U}_k}$ is a proper submersion to $\widehat{U}_{k-1}\times\mathbb R_+.$

After composing these maps $n$ times, one gets the following Lemma:

\begin{lem}
$(\rho_{\breve V_0},\dots,\rho_{\breve V_{k-1}}):\widehat{U}_k\to(\mathbb R_+)^k$ is a proper
submersion. Therefore, $\widehat{U}_k$ is diffeomorphic to $\tau_k\times (\mathbb R_+^k).$
\end{lem}

Consider the preimages $U_k=(p_n\circ\dots\circ p_k)^{-1}(\widehat{U}_k)\subset\widetilde{W}_k.$
By Lemma \ref{Stratification of the normalization diagram}, $U_k\subset
\breve{\widetilde{W}}_k$ and $(p_n\circ\dots\circ p_k)|_{U_k}:U_k\to\widehat{U}_k$ is a covering.

Denote $\overline{U}_k:=U_k\cup p_{k-1}(U_{k-1})$ for $k=n,n-1,\dots,1.$

\begin{lem}\label{open subset in the set of regular points}
$\overline{U}_k\subset\widetilde{W}_k$ is an open subset consisting of regular points of
$\widetilde{W}_k.$

\end{lem}

\begin{proof}
$\widehat{U}_k\cup \widehat{U}_{k-1}$ is an open subset in ${\breve V_k}\cup {\breve V_{k-1}}.$ Indeed, it is the preimage of the $\widehat{U}_{k-1}$ under the restriction of the projection $\pi_{\breve{V_{k-1}}},$ restricted to $U_{\breve V_{k-1}}\cap (\breve V_{k-1}\cup\breve V_k).$

In turn, $\overline{U}_k$ is the preimage of $\widehat{U}_k\cup \widehat{U}_{k-1}$ under $(p_n\circ\dots\circ p_k)|_{\breve{\widetilde{W}}_k \cup
{\breve W_{k-1}}}.$

Also, by Lemma \ref{(n-1)-dimensional stratum}, ${\breve W_{k-1}}$ consists of
regular points of $\widetilde{W}_k.$
\end{proof}

We need the following Lemma about lifting the control data:

\begin{lem}\label{lifting control data}
Let $V$ and $V'$ be two stratified spaces, consisting of two strata each: $V=X\sqcup Y,\ $ $X<Y,$
and $V'=X'\sqcup Y',\ $ $X'<Y'.$ Let $p: V'\to V$ be a map, such that $p|_{X'}$ is a covering
over $X$ and $p|_{Y'}$ is a covering over $Y.$ Let $U_X\subset V,\ \pi_X:U_X\to X,$ and
$\rho_X:U_X\to\mathbb R_{\ge 0}$ be the control data on $V.$ Then there exist control data
$U_{X'},\pi_{X'},\rho_{X'}$ on $V',$ such that
\begin{enumerate}
\item $\rho_X\circ p=\rho_{X'};$
\item $\pi_X\circ p=p\circ \pi_{X'}.$
\end{enumerate}
\end{lem}

\begin{proof}
We set the tubular neighborhood $U_{X'}:=p^{-1}(U_X).$ The tubular function $\rho_{X'}$ is defined
by the property (1). The retraction $\rho_{X'}$ is defined uniquely by the property (2) and
continuity.
\end{proof}

We apply the Lemma \ref{lifting control data} to the $V=\widehat{U}_k\sqcup \widehat{U}_{k-1}$ and
$V'=\overline{U}_k=p_{k-1}(U_{k-1})\sqcup U_k.$ Let $\pi_{p_{k-1}(U_{k-1})}:\overline{U}_k\to
p_{k-1}(U_{k-1})$ and $\rho_{p_{k-1}(U_{k-1})}:\overline{U}_k\to\mathbb R_{\ge 0}$ be the corresponding
retraction and tubular function. We have the following corollary:

\begin{cor}
For any $k=n,n-1,\dots,1$ the connected components of $U_k$ are in natural one-to-one
correspondence with the connected components of $U_{k-1}.$
\end{cor}

\begin{proof}
Indeed, the map from the connected components of $U_k$ to the connected components of
$p_{k-1}(U_{k-1})$ is given by the retraction $\rho_{p_{k-1}(U_{k-1})}.$ Existence of the
inverse to this map follows from the fact that $p_{k-1}(U_{k-1})\subset\overline{U}_k$ is a complex
hypersurface in the manifold $\overline{U}_k.$ Indeed, if $H\subset M$ is a hypersurface in a complex manifold $M,$ then there is only one connected component of $M$ in a neighborhood of a connected component of $H.$

Finally, $p_{k-1}|_{U_{k-1}}$ is an isomorphism to the image.
\end{proof}

Pick a point $a_{\alpha}\in W_0.$ Let $U_1^{\alpha},\dots,U_n^{\alpha}$ be the corresponding
connected components of $U_1,\dots,U_n$ correspondingly. Let also
$\overline{U}_k^{\alpha}:=U_k^{\alpha}\cup p_{k-1}(U_{k-1}^{\alpha})$ be the corresponding connected
components of $\overline{U}_k.$

Let $\widetilde{\tau}_k:=(p_n\circ\dots\circ p_k)^{-1}(\tau_k).$ Note, that
$\widetilde{\tau}_k\subset U_k$ is a union of connected components, one in each $U_k^{\alpha}.$ Let $\widetilde{\tau}_k^{\alpha}\subset U_k^{\alpha}$ be the corresponding connected component.

\begin{lem}\label{induction step}
$\phi_{p_{k-1}(U_{k-1}),U_k}\circ
(p_{k-1}|_{U_{k-1}})_*([\widetilde{\tau}_{k-1}^{\alpha}])=[\widetilde{\tau}_k^{\alpha}].$
\end{lem}

\begin{proof}
Since we have chosen the control data on $\overline{U}_k$ to be coherent with the control data on $\widehat{U}_k\cup \widehat{U}_{k-1}$ (given by restricting the control data from the ambient space), we have the equality on the level of representatives.
\end{proof}

Now we use the local parameters $(u_1,\dots,u_n)$ to construct cycles $\gamma_k^{\alpha}\subset
U_k^{\alpha}$ such that, on the one side, it is obvious that
$$
res_{\{F,a_{\alpha}\}}(\omega)=\frac{1}{(2\pi i)^n}\int\limits_{\gamma_n^{\alpha}}\omega,
$$
\noindent and, on the other, $\gamma_k^{\alpha}$ is homologically equivalent to
$\widetilde{\tau}_k^{\alpha}$ in $U_k^{\alpha}.$

Function $u_k$ is regular and non-vanishing in $U_k^{\alpha}\subset\breve{\widetilde{W}_k}$ and has
zero of order one at a generic point of $p_{k-1}(U_{k-1})\subset\overline{U}_k.$ It follows
immediately, that $u_k$ is regular on $\overline{U}_k$ and the equation $u_k=0$ defines
$p_{k-1}(U_{k-1})$ in $\overline{U}_k.$

The following Lemma easily follows from the above observation:

\begin{lem}
There exist smooth positive real functions $\epsilon_1,\dots, \epsilon_n,\ $ $\epsilon_k:\mathbb
C^{k-1}\to \mathbb R_+,$ and open subsets $B_k\subset U_k,\ $ $k=1,\dots,n,$ such that
$$
(u_1,\dots,u_k): B_k\to A_k:=\{(z_1,\dots,z_k): |z_i|<\epsilon_i(z_1,\dots,z_{i-1}),\
i=1,\dots,k\}\subset\mathbb C^k
$$
\noindent are biholomorphisms. (Note, that $\epsilon_1$ is a constant.)
\end{lem}

Let $\delta_1,\dots, \delta_n\in\mathbb R_+$ be small enough, so that $\{(z_1,\dots,z_n):
|z_i|=\delta_i,\ i=1,\dots,n\}\subset A_n.$

\begin{df}
Denote $\gamma_k^{\alpha}=\{x\in B_k: |u_i(x)|=\delta_i,\ i=1,\dots,k\},$ and
$\gamma_0^{\alpha}=a_{\alpha}.$
\end{df}

It follows immediately from the definition of the Parshin residue, that
$$
res_{\{F,a_{\alpha}\}}(\omega)=\frac{1}{(2\pi i)^n}\int\limits_{\gamma_n^{\alpha}}\omega.
$$

\begin{lem}
$\gamma_k^{\alpha}$ and $\widetilde{\tau}_k^{\alpha}$ define the same homology class in
$H_k(U_k).$
\end{lem}

\begin{proof}
We prove this Lemma by induction. For $k=0$ one has
$\gamma_0^{\alpha}=\widetilde{\tau}_0^{\alpha}=a_{\alpha}.$

For the induction step, one uses the Lemma \ref{induction step} and the similar observation for
cycles $\gamma_k^{\alpha}.$
\end{proof}

So, we proved the following Theorem:
\begin{te}\label{Residue via coboundary homomorphisms}
$$
res_F(\omega)=\frac{1}{(2\pi i)^n}\int\limits_{\Delta_F}\omega,
$$
\noindent where $F:=\{V_n\supset\dots\supset V_0\}$ and $\Delta_F=\phi_{\breve V_{n-1},\breve
V_n}\circ\dots\circ\phi_{\breve V_0,\breve V_1}([V_0])\in H_n(\breve V_n).$
\end{te}

In order to get the Parshin's Reciprocity Law from the Theorem \ref{Relation} and the above
consideration, one needs a fixed stratification of $V$, such that all non-zero residues of a given
form $\omega$ are in the flags consisting of closures of strata of the stratification. It turns
out, that any Whitney stratification, such that $\omega$ is regular on the top dimensional stratum
is good enough. More precisely, we have the following Theorem:

\begin{te}\label{criteria for non-trivial residues}
Let $V$ be an $n$-dimensional variety and $\omega$ be a meromorphic $n$-form on $V.$ Let ${\bf
S_{\omega}}$ be a Whitney stratification of $V,$ such that $\omega$ is regular on ${\breve V}$.

Let $F=\{V_n\supset\dots\supset V_0\}$ be a flag of irreducible subvarieties of $V,$ $\dim V_i=i.$
Suppose that at least one of $V_i$'s is not the closure of a stratum of ${\bf S_{\omega}}.$ Then
$res_{\{F,a_{\alpha}\}}\omega=0$ for all $a_{\alpha}\in W_0.$
\end{te}

\begin{proof}
Consider the normalization diagram for the flag $F.$ Let $a_{\alpha}\in W_0,$ and let $(u_1,\dots,u_n)$
be local parameters. Let ${\bf S}$ be a stratification of $V$ satisfying conditions of the Lemma
\ref{Stratification of the normalization diagram}, and such that all strata of the stratification
${\bf S_{\omega}}$ are unions of strata of ${\bf S}.$ As usual, we denote by ${\breve V_k}$ the
stratum of ${\bf S}$ which is open and dense in $V_k.$

The proof of the theorem is based on two observations:

\begin{enumerate}
\item Let $X'<Y'$ be consecutive strata of ${\bf S_{\omega}},$ and let $X<Y$ be the
consecutive strata of ${\bf S},$ such that $X$ is an open dense subset in $X',$ and $Y$ is an open
dense subset in $Y'.$ Let $i_X:X\hookrightarrow X'$ and $i_Y:Y\hookrightarrow Y'$ be the
embeddings. Then it easily follows from the construction of coboundary operators $\phi_{X,Y}$ and
$\phi_{X',Y'}$ and the independence of these operators from the choice of the control data, that
$\phi_{X',Y'}\circ i_{X*}=i_{Y*}\circ \phi_{X,Y}.$
\item Let $Y'$ be a stratum of ${\bf S_{\omega}}$ and let $X<Y$ be consecutive strata of ${\bf S},$
such that $(X\cup Y)\subset Y'$ and $Y$ is open and dense in $Y'.$ Then $i_*\circ \phi_{X,Y}=0,$
where $i:Y\hookrightarrow Y'$ is the embedding. Moreover, if $A$ is a representative of a homology
class in $H_*(X)$ and $B$ is the representative of the $\phi_{X,Y}([A])\in H_{*+\dim Y-\dim
X-1}(Y),$ constructed in the standard way, then every connected component of $B$ is homologically
equivalent to $0$ in $Y.$ Indeed, one can use the control data  on ${\bf S}$ to embed the mapping
cone of $\pi_{X}|_B:B\to A$ into $Y'.$
\end{enumerate}

Let $k$ be the largest number, such that ${\breve V_k}$ is a subset of a stratum of ${\bf
S_{\omega}}$ of dimension bigger then $k.$ For $m=k+1,\dots,n$, let ${\breve V_m}'$ be the stratum
of ${\bf S_{\omega}},$ such that ${\breve V_m}\subset {\breve V_m}'.$ Note, that $\dim {\breve
V_m}'=m$ and ${\breve V_m}$ is open and dense in ${\breve V_m}'.$ Moreover, by dimension reasons
and the condition of the frontier, ${\breve V_k}\subset {\breve V_{k+1}}'.$

Let $i_m:{\breve V_m}\hookrightarrow {\breve V_m}'$ be the embedding. Then, according to the first
observation, one has

$$
(i_n)_*\circ\phi_{{\breve V_{n-1}},{\breve V_n}}\circ\dots\circ \phi_{{\breve V_1},{\breve V_0}}=
\phi_{{\breve V_{n-1}}',{\breve V_n}'}\circ\dots\circ \phi_{{\breve V_{k+1}}',{\breve
V_{k+2}}'}\circ (i_{k+1})_*\circ \phi_{{\breve V_k},{\breve V_{k+1}}}\circ\dots\circ \phi_{{\breve
V_1}{\breve V_0}}.
$$

On the other side, according to the second observation, $(i_{k+1})_*\circ \phi_{{\breve V_k},{\breve
V_{k+1}}}=0.$ Therefore,

$$
(i_n)_*\circ\phi_{{\breve V_{n-1}},{\breve V_n}}\circ\dots\circ \phi_{{\breve V_1}{\breve V_0}}=0,
$$

\noindent and, since $\omega$ is regular in ${\breve V_n}',$ $res_{F} \omega=0.$ Moreover, easy to
see, that every connected component of the standard representative of the $\phi_{{\breve
V_{n-1}},{\breve V_n}}\circ\dots\circ \phi_{{\breve V_1}{\breve V_0}}([V_0])$ is homologically
equivalent to $0.$ Therefore, $res_{F,a} \omega=0$ for any $a\in W_0.$
\end{proof}

\noindent{\bf Corollary.} There are only finitely many non-zero Parshin's residues for a given
meromorphic form.

Note, that the Parshin Reciprocity Law follows from the Theorems \ref{Residue via coboundary
homomorphisms}, \ref{criteria for non-trivial residues}, and \ref{Relation}.

\bibliographystyle{plain}

\end{document}